\documentclass{elsart3-1-ppt}

\usepackage{amssymb}

\usepackage[english,francais]{babel}

\newtheorem{theorem}{Theorem}[section]

\newtheorem{e-proposition}[theorem]{Proposition}

\newtheorem{e-definition}[theorem]{Definition\rm}

\newtheorem{theoreme}{Th\'eor\`eme}

\newtheorem{proposition}[theoreme]{Proposition}

\renewcommand{\theequation}{\arabic{equation}}
\newcommand{\R}{\mathbb{R}}
\newcommand{\eps}{\varepsilon}

\newcommand{\nx}{\nabla_x}
\newcommand{\nv}{\nabla_v}

\newcommand{\rhoF}{\rho_F}

\newcommand{\mF}{m_F}
\newcommand{\href}[1]{(\ref{#1})}
\renewcommand{\a}{\mathsf{a}}
\renewcommand{\b}{\mathsf{b}}
\newcommand{\A}{\mathsf{A}}
\newcommand{\T}{\mathsf{T}}
\renewcommand{\L}{\mathsf{L}}

\newcommand{\nrm}[1]{\|#1\|}

\newcommand{\eqref}[1]{(\ref{#1})}
\newcommand{\be}[1]{\begin{equation}\label{#1}}
\newcommand{\ee}{\end{equation}}
\newenvironment{multline*}{\begin{eqnarray*}}{\end{eqnarray*}}
\newenvironment{equation*}{\[}{\]}

\def\og{\leavevmode\raise.3ex\hbox{$\scriptscriptstyle\langle\!\langle$~}}
\def\fg{\leavevmode\raise.3ex\hbox{~$\!\scriptscriptstyle\,\rangle\!\rangle$}}

\journal{the Acad\'emie des sciences}
\begin{document}

\centerline{}
\begin{frontmatter}

\selectlanguage{english}


\title{Hypocoercivity for kinetic equations with linear relaxation
  terms}

\selectlanguage{english}
\author[Ceremade]{Jean Dolbeault}
\ead{dolbeaul@ceremade.dauphine.fr}
\author[Ceremade]{Cl\'ement Mouhot}
\ead{Clement.Mouhot@ceremade.dauphine.fr}
\author[Schmeiser]{Christian Schmeiser}
\ead{Christian.Schmeiser@univie.ac.at}

\address[Ceremade]{Ceremade (UMR CNRS no. 7534), Universit\'e Paris-Dauphine, Place de Lattre de Tassigny, 75775 Paris 16, France}
\address[Schmeiser]{Fakult\"at f\"ur Mathematik, Universit\"at Wien,
  Nordbergstra{\ss}e 15, 1090 Wien, Austria}
  
\begin{abstract}
  \selectlanguage{english} This note is devoted to a simple method for
  proving the hypocoercivity associated to a kinetic equation
  involving a linear time relaxation operator. It is based on the
  construction of an adapted Lyapunov functional satisfying a
  Gronwall-type inequality. The method clearly distinguishes the
  coercivity at microscopic level, which directly arises from the
  properties of the relaxation operator, and a spectral gap inequality
  at the macroscopic level for the spatial density, which is connected
  to the diffusion limit. It improves on previously known results. Our
  approach is illustrated by the linear BGK model and a relaxation
  operator which corresponds at macroscopic level to the linearized
  fast diffusion.

\vskip 0.5\baselineskip

\selectlanguage{francais}
\noindent{\bf R\'esum\'e} \vskip 0.5\baselineskip \noindent {\bf
  Hypocoercivit\'e pour des \'equations cin\'etiques avec termes de
  relaxation lin\'eaires}

Cette note est consacr\'ee \`a une m\'ethode simple pour d\'emontrer
l'hypocoercivit\'e associ\'ee \`a une \'equation cin\'etique
contenant un op\'erateur de relaxation lin\'eaire; il s'agit de
construire une fonctionnelle de Lyapunov adapt\'ee v\'erifiant une
in\'egalit\'e de type Gronwall. La m\'ethode distingue clairement la
coercivit\'e au niveau microscopique, qui provient directement des
propri\'et\'es de l'op\'erateur de relaxation, et une in\'egalit\'e de
trou spectral pour la densit\'e spatiale, qui est reli\'ee \`a la limite
de diffusion. Elle am\'eliore les r\'esultats ant\'erieurs.  Notre
approche est illustr\'ee par le mod\`ele de BGK lin\'eaire et par un
op\'erateur de relaxation qui correspond, au niveau macroscopique, \`a
la diffusion rapide lin\'earis\'ee.

\end{abstract}

\end{frontmatter}

\selectlanguage{francais}
\section*{Version fran\c{c}aise abr\'eg\'ee}

Cette note est une contribution \`a la th\'eorie de
l'\emph{hypocoercivit\'e}, voir
\cite{MR2339441,Mem-villani,Mouhot-Neumann}, dont le but est d'estimer
des taux exponentiels de retour \`a un \'equilibre global pour des
\'equations cin\'etiques dans lesquelles le terme de collision ne
contr\^ole que le retour \`a un \'equilibre local. Cette question a
\'et\'e partiellement r\'esolue pour l'\'equation, hypoelliptique, de
Vlasov-Fokker-Planck par F. H\'erau and F. Nier dans \cite{Herau-Nier}, puis par C.~Villani dans \cite{Mem-villani}, et pour le mod\`ele de BGK
lin\'eaire, non hypoelliptique, dans \cite{Herau} par F. H\'erau. 
En introduisant de nouveaux op\'erateurs, nous am\'eliorons et nous
simplifions les r\'esultats ant\'erieurs, tout en mettant en
\'evidence la relaxation \`a l'\'echelle microscopique due au terme de
collision, ici un op\'erateur de relaxation lin\'eaire, et la
relaxation \`a l'\'echelle macroscopique, qui vient d'une
propri\'et\'e de trou spectral (voir \cite{BlBoDoGrVa}) d'un
op\'erateur li\'e \`a la limite de diffusion (voir \cite{DoMaOeSc}) et
portant sur les densit\'es spatiales uniquement.

\medskip Soit $V$ un \emph{potentiel ext\'erieur} sur $\R^d$, $d\ge1$.
Consid\'erons l'\'equation cin\'etique
\begin{equation}
\label{eq:baseF}
\partial_t f + \T\,f = \L\,f\,, \quad f = f(t,x,v)\,, \ x \in \R^d, \ v \in \R^d,
\end{equation}
o\`u $\T := v \cdot \nx - \nx V \cdot \nv$ est un op\'erateur de
transport et o\`u l'op\'erateur de relaxation lin\'eaire $\L$ est
d\'efini par
\[
\L\,f=\Pi\,f-f\,,\quad\Pi\,f := \frac\rho\rhoF\,F(x,v)\,, \ F(x,v) >0,
\quad \rho=\rho(f) := \int_{\R^d} f \, dv\quad\mbox{et}\quad
\rhoF=\rho(F)\,.
\]
Sur $\R^d\times\R^d\ni(x,v)$, on notera
$d\mu(x,v)=F(x,v)^{-1}\,dx\,dv$ o\`u $F$ est une mesure de probabilit\'e strictement positive. Norme
et produit scalaire seront par d\'efaut ceux de $L^2(d\mu)$. La
donn\'ee initiale $f_0\in L^2(d\mu)$ est normalis\'ee par
$\int_{\R^d\times\R^d}f_0\,dx\,dv=1$, et la solution du probl\`eme de
Cauchy est unique car $T$ est antisym\'etrique et $L$ est sym\'etrique
n\'egatif.

\medskip Le premier cas d'application de notre m\'ethode est le cas
d'un \'equilibre Maxwellien: $F(x,v):=M(v)\,e^{-V(x)}$ avec
$M(v):=(2\pi)^{-d/2} \,e^{-|v|^2/2}$.\smallskip
\begin{theoreme}\label{Thm:Main1F} 
  Supposons que $V\in W^{2,\infty}_{\rm loc}(\R^d)$, $d\ge 1$, est tel
  que: 1) $\int_{\R^d}e^{-V}dx=1$, 2) il existe une constante
  $\Lambda>0$ telle que
  $\int_{\R^d}|u|^2\,e^{-V}dx\leq\Lambda\int_{\R^d}|\nx
  u|^2\,e^{-V}dx$ pour tout $u\in H^1(e^{-V}dx)$ v\'erifiant
  $\int_{\R^d}u\,e^{-V}dx=0$, 3) il existe des constantes $c_0>0$,
  $c_1>0$ et $\theta\in (0,1)$ telles que $\Delta V\leq
  \frac\theta2\,|\nx V(x)|^2+c_0$ et $|\nx^2 V(x)|\leq c_1\,(1+|\nx
  V(x)|)\;\forall\,x\in\R^d$, 4) $\int_{\R^d}|\nx
  V|^2\,e^{-V}dx<\infty$.

Pour tout $\eta>0$, il existe une constante positive $\lambda=\lambda(\eta)$, explicite, pour
  laquelle toute solution de \eqref{eq:baseF} dans $L^2(d\mu)$
  v\'erifie:
  \[
  \forall\,t\ge 0\,,\quad\|f(t)-F\|^2\le(1+ \eta)\,\|f_0-F\|^2\,e^{-\lambda t}.
  \]
\end{theoreme}

Le deuxi\`eme cas d'application correspond \`a un mod\`ele reli\'e \`a
l'\'equation de diffusion rapide avec terme de d\'erive,
lin\'earis\'ee autour de profils stationnaires de type Barenblatt. On
supposera que
\[
F(x,v):=\omega\left(\frac 12\,|v|^2+V(x)\right)^{-(k+1)},\quad V(x)=\left(1+|x|^2\right)^\beta
\]
pour simplifier la pr\'esentation. Ici, $\rhoF=\omega_0\,V^{d/2-k-1}$,
$\omega$ et $\omega_0$ sont des constantes de normalisation,
positives. Des choix plus g\'en\'eraux de $V$ ainsi que les
hypoth\`eses correspondantes seront donn\'es dans \cite{DMS-1}.
\smallskip
\begin{theoreme}\label{Thm:Main2F} Soit $d\ge 1$, $k>d/2+1$. Il existe une constante $\beta_0>1$
  telle que, pour tout $\beta\in(\min\{1,(d-4)/(2k-d-2)\},\beta_0 )$, il existe deux constantes strictement positives $C$ et 
  $\lambda$, explicites, pour laquelle toute solution de
  \eqref{eq:baseF} dans $L^2(d\mu)$ v\'erifie:
  \[
  \forall\,t\ge 0\,,\quad\|f(t)-F\|^2\le C\,\|f_0-F\|^2\,e^{-\lambda t}.
  \]\end{theoreme}

Sur $L^2(d\mu)$, \`a l'aide de l'op\'erateur de projection $\Pi$, on d\'efinit les op\'erateurs $
  \b f := \Pi\, ( v \, f)$, $\a f := \b\, ( \T\,f)$, $\hat \a\, f := - \Pi\,(\nx f)$, $\A := (1 + \hat \a\cdot \a\,\Pi)^{-1}\, \hat \a\cdot \b$
et la fonctionnelle
\[
H(f) := \frac{1}{2}\,\|f\|^2 + \eps\,\langle \A\,f,f\rangle\,.
\]
Si $f$ est une solution de l'\'Eq. \eqref{eq:baseF}, alors $D(f-F):=\frac
d{dt}H(f-F)$ est donn\'e par
\[
D(f):=\langle f,\L\,f\rangle-\eps\,\langle
\A\,\T\,\Pi\,f,f\rangle-\eps\,\langle
\A\,\T\,(1-\Pi)\,f,f\rangle+\eps\,\langle
\T\,\A\,f,f\rangle+\eps\,\langle \L\,f,(\A+\A^*)\,f\rangle\,.
\]
Pour simplifier les notations, on remplace $f$ par $f-F$ de sorte que $\int_{\R^d}f\,dx\,dv=0$. La preuve des Th\'eor\`emes~\ref{Thm:Main1F} et \ref{Thm:Main2F} consiste alors \`a montrer que $D(f)+\lambda\,H(f)\ge 0$. Elle repose
principalement sur deux estimations. D'une part, par l'in\'egalit\'e
de trou spectral, le terme
\[
-\eps\, \langle
\A\,\T\,\Pi\,f,f\rangle\le-\eps\,\frac\Lambda{1+\Lambda}\,\|\Pi\,f\|^2
\]
contr\^ole les termes macroscopiques \`a l'ordre $\eps$,
c'est-\`a-dire proportionnels \`a $\|\Pi\,f\|^2$. D'autre part, on
remarque que
\[
2\,\|\A\,f\|^2+\|\T\,\A\,f\|^2\le\|(1-\Pi)\,f\|^2.
\]
Le terme $\langle f,\L\,f\rangle=-\|(1-\Pi)\,f\|^2$ permet alors de
contr\^oler tous les autres termes, et en particulier
$\|(\A\,\T\,(1-\Pi))f\|^2$ que l'on \'evalue en consid\'erant
l'op\'erateur dual: si
$(\A\,\T\,(1-\Pi))^*f=(\hat\a\cdot\a\,(1-\Pi))^*g$ avec $g=\left( 1 +
  \hat\a \cdot\a\, \Pi \right)^{-1}f$, alors $u:=\rho(g)/\rhoF$ est
donn\'ee comme solution de l'\'equation \be{Eqn:u1F}
\rho(f)=\rhoF\,u-\frac 1d\,\nx\left(\mF\,\nx u\right) \ee avec
$\mF(x):=\int_{\R^d}|v|^2\,F(x,v)\,dv$, et il suffit alors d'\'evaluer $\|(\A\,\T\,(1-\Pi))^*f\|^2$ gr\^ace \`a \eqref{Eqn:Dual}. Ceci revient essentiellement \`a \'etablir une estimation $H^2$ pour
la solution de \eqref{Eqn:u1F}. Dans le cas Maxwellien, il faut
d'abord \'etablir une in\'egalit\'e de Poincar\'e am\'elior\'ee: il
existe une constante $\kappa>0$ telle que, pour tout $u\in
H^1(e^{-V}dx)$ v\'erifiant $\int_{\R^d}u\,e^{-V}dx=0$,
$\kappa\,\int_{\R^d}|\nx V|^2\,|u|^2\,dx\le\|\nx u\|_0^2$, d'o\`u l'on
d\'eduit d'abord que $\int_{\R^d}|\nx V|^2\,|\nx u|^2\,e^{-V}dx$ est
born\'e par $\|f\|^2$, puis que $\nrm{\nx^2u}_0^2$ est aussi
contr\^ol\'e par $\|f\|^2$. Dans le cas de la diffusion rapide, il
suffit de multiplier \eqref{Eqn:u1F} par $V^{1-1/\beta}u$ et par $V\,\Delta u$, 
puis d'effectuer quelques int\'egrations par parties, pour contr\^oler en
d\'efinitive $\|(\A\,\T\,(1-\Pi))^*f\|^2$ par $\|f\|^2$, ce qui permet
de conclure.

\selectlanguage{english}
\setcounter{equation}{0}
\par\medskip\centerline{\rule{2cm}{0.2mm}}\medskip
\section{Introduction}\label{Sec:Intro}

A fundamental question which goes back to the early days of kinetic
theory is to estimate the \emph{rate of relaxation} of the solutions
towards a global equilibrium. This is not an easy issue since the
collision term responsible for the relaxation acts, in most of the
cases, only on the velocity space. Rates of convergence have been
investigated in many papers for the so-called homogeneous kinetic
equations, but understanding how the transport operator interacts with
collisions to produce a global relaxation is a different and much more
recent story. The point is to understand how the spatial density evolves towards a density corresponding to a
distribution function which is simultaneously in the kernels of the
collision and transport operators, a property of the stationary
solutions of many kinetic equations. There is an obvious link with
diffusion or hydrodynamic limits. A key feature of our approach is
that it clearly distinguishes the mechanisms of relaxation at
\emph{microscopic level} (convergence towards a local equilibrium, in
velocity space) and \emph{macroscopic level} (convergence of the
spatial density to a steady state), where the rate is given by a
spectral gap which has to do with the underlying diffusion equation
for the spatial density. See \cite{DMS-1} for more details.

First non constructive results were obtained by Ukai \emph{et al.,\/} see for instance \cite{Ukai}. Constructive methods inspired
from \emph{hypoelliptic theory} (see {\it e.g.}  \cite{Hormander})
were then brought into the field of kinetic theory by F. H\'erau and
F.  Nier, see for instance \cite{Herau-Nier} in case of the
Vlasov-Fokker-Planck equation. In a recent paper, \cite{Herau}, F.
H\'erau studied with such tools the case of an operator of zeroth
order in the derivatives, which is known in the kinetic literature as
the linear Boltzmann relaxation operator. Our approach is done in the
spirit of \cite{Herau} but in a simplified framework for which the
order of the operator plays no role. Moreover, explicit estimates on
the relaxation rate easily follow, weaker assumption on the external
potential than in \cite{Herau} are needed and the method applies to
more general relaxation operators of which we shall give an example.
This example is based on kinetic equations which have been studied in
\cite{DoMaOeSc} and give equations of fast diffusion in the diffusion
limit.

The hypoelliptic theory is mainly focused on the regularization
properties of the evolution operator, but in some cases the hypoelliptic estimates also imply a result of relaxation
to equilibrium. However both questions are independent and have to be
distinguished. In the \emph{hypocoercivity} approach, the purpose is
centered on the asymptotic behavior and the quantification of the
relaxation rates. More precisely our goal is to construct a Lyapunov
functional, or generalized entropy, and establish an inequality
relating the entropy and its time derivative along the flow defined by
the evolution equation.  To establish the inequality is then
equivalent to prove an exponential rate.  Such an approach has
systematically been tackled by C. Villani, see
\cite{MR2339441,Mem-villani}, and has been successfully applied to
various models, see for instance \cite{Mouhot-Neumann}. It is also
related to recent works on non-linear Boltzmann and Landau equations,
see {\it e.g.} \cite{Guo}. In this note, we develop a new approach based on
operators with less algebraic properties than the ones of the
hypoelliptic theory, but which are better adapted to the micro-macro
decomposition of the distribution function and give a much simpler
insight of the mechanisms responsible of the relaxation at both
levels. We illustrate our approach on two
examples: the linear BGK and the linearized fast diffusion models. We
refer the interested reader to a forthcoming paper, \cite{DMS-1}, in
which the theory will be developed at a more general and abstract
level.

\section{Main results}\label{Sec:Main}

Let $V$ be a given \emph{external potential} on $\R^d$, $d\ge1$, and
consider the kinetic equation
\begin{equation}
\label{eq:base}
\partial_t f + \T\,f = \L\,f\,, \quad f = f(t,x,v)\,, \ x \in \R^d, \ v \in \R^d,
\end{equation}
where $\T := v \cdot \nx - \nx V \cdot \nv$ is a transport operator,
and the linear relaxation operator $\L$ is defined by
\[
\L\,f=\Pi\,f-f\,,\quad\Pi\,f := \frac\rho\rhoF\,F(x,v)\,, \quad
\rho=\rho(f) := \int_{\R^d} f \, dv\quad\mbox{and}\quad
\rhoF=\rho(F)
\]
for some function $F(x,v)>0$ which only depends on $|v|^2/2 + V(x)$. On
$\R^d\times\R^d\ni(x,v)$, we consider the measure
$d\mu(x,v)=F(x,v)^{-1}\,dx\,dv$ where $F$ is a positive probability measure. Unless it is
explicitly specified, the scalar product and the norm are the ones of
$L^2(d\mu)$: $\langle
f,g\rangle=\int_{\R^d\times\R^d}f\,g\,d\mu$ and
$\|f\|^2=\langle f,f\rangle$. Throughout this paper, Eq.
\eqref{eq:base} is supplemented with a nonnegative initial datum
$f_0\in L^2(d\mu)$ such that $\int_{\R^d\times\R^d}f_0\,dx\,dv=1$. We
shall assume that a unique solution globally exists. This is granted
under additional technical assumptions, see for instance
\cite{DoMaOeSc}. The goal of this note is to state hypocoercivity
results in the two following cases.

\subsection{Maxwellian case}
We assume that $F(x,v):=M(v)\,e^{-V(x)}$ with
$M(v):=(2\pi)^{-d/2}\,e^{-|v|^2/2}$, where $V(x)=C_k\,(1+|x|^2)^{k/2}$ for $k>1$ and $C_k$ is appropriately chosen or, more generally, satisfies the following assumptions:\begin{itemize}
\item[(H1)] {\sl Regularity:\/} $V\in W^{2,\infty}_{\rm loc}(\R^d)$.
\item[(H2)] {\sl Normalization:\/} $\int_{\R^d}e^{-V}dx=1$.
\item[(H3)] {\sl Spectral gap condition:\/} there exists a positive
  constant $\Lambda$ such that
  $\int_{\R^d}|u|^2\,e^{-V}dx\leq\Lambda\int_{\R^d}|\nx
  u|^2\,e^{-V}dx$ for any $u\in H^1(e^{-V}dx)$ such that
  $\int_{\R^d}u\,e^{-V}dx=0$.
\item[(H4)] {\sl Pointwise condition 1:\/} there exists $c_0>0$ and
  $\theta\in (0,1)$ such that $\Delta V\leq \frac\theta2\,|\nx
  V(x)|^2+c_0\;\forall\,x\in\R^d$.
\item[(H5)] {\sl Pointwise condition 2:\/} there exists $c_1>0$ such
  that $|\nx^2 V(x)|\leq c_1\,(1+|\nx V(x)|)\;\forall\,x\in\R^d$.
\item[(H6)] {\sl Growth condition:\/} $\int_{\R^d}|\nx
  V|^2\,e^{-V}dx<\infty$.
\end{itemize}\smallskip

\begin{theorem}\label{Thm:Main1} For any $\eta>0$, there exists an explicit, positive
  constant $\lambda=\lambda(\eta)$ such that, under the above assumptions, the
  solution of~\eqref{eq:base} satisfies:
  \[
  \forall\,t\ge 0\,,\quad\|f(t)-F\|^2\le(1+ \eta)\,\|f_0-F\|^2\,e^{-\lambda t}.
  \]
\end{theorem}

\subsection{Fast diffusion case}
For some $\beta>0$ to be specified later, we assume that
\[
F(x,v):=\omega\left(\frac 12\,|v|^2+V(x)\right)^{-(k+1)},\quad
V(x)=\left(1+|x|^2\right)^\beta
\]
where $\omega$ is a normalization constant chosen such that
$\int_{\R^d\times\R^d} F\,dx\,dv=1$ and $\rhoF=\omega_0\,V^{d/2-k-1}$
for some $\omega_0>0$. More general choices for~$V$ and corresponding
assumptions can be found in \cite{DMS-1}. \smallskip
\begin{theorem}\label{Thm:Main2} Let $d\ge 1$, $k>d/2+1$. There exists 
  a constant $\beta_0>1$ such that, for any $\beta\in(\min\{1,(d-4)/(2k-d-2)\},\beta_0)$, there are two positive, explicit
  constants $C$ and $\lambda$ for which the solution of \eqref{eq:base}
  satisfies:
  \[
  \forall\,t\ge 0\,,\quad\|f(t)-F\|^2\le C\,\|f_0-F\|^2\,e^{-\lambda t}.
  \]
\end{theorem}

\subsection{A Lyapunov functional}
On $L^2(d\mu)$, $\Pi$ is the orthogonal projection onto the space of local equilibria. Let us define
\begin{equation*}
  \b f := \Pi\, ( v \, f)\,,\quad\a f := \b\, ( \T\,f)\,,\quad\hat \a\, f := - \Pi\,(\nx f)\quad\mbox{and}\quad\A := (1 + \hat \a\cdot \a\,\Pi)^{-1}\, \hat \a\cdot \b\,.
\end{equation*}
In the definition of $\A$, we take the product coordinate by coordinate. These operators can be rewritten as
\[
\b\,f=\frac F\rhoF\,\int_{\R^d}v\,f\,dv\,,\quad\a\,f=\frac
F\rhoF\,\left(\nx\cdot\int_{\R^d}v\otimes v\,f\,dv+\rho(f)\,\nx
  V\right),\quad\mbox{and}\quad\hat\a\,f=-\frac F\rhoF\,\nx\rho(f)\,.
\]
Let us also note that $\A\,\T=(1 + \hat \a\cdot \a\,\Pi)^{-1}\, \hat
\a\cdot \a$ and $\a\,\Pi\,f=\frac
F\rhoF\,\frac\mF{d}\,\nx\!\left(\frac{\rho(f)}\rhoF\right)$ where
$\mF:=\int_{\R^d} |v|^2 \, F(\cdot,v) \, dv$ $=d\int_{\R^d} |v_i|^2 \,
F(\cdot,v) \, dv$ for any $i=1$, $2$\ldots $d$. Define the functional
\[
H(f) := \frac{1}{2}\,\|f\|^2 + \eps\,\langle \A\,f,f\rangle\,.
\]
The operator $\T$ is skew-symmetric on $L^2(d\mu)$. If $f$ is a
solution of Eq. \eqref{eq:base}, then
\begin{eqnarray*}
  &&\frac d{dt}H(f-F)=D(f-F)\\
  &&\mbox{with}\quad D(f):=\langle f,\L\,f\rangle-\eps\,\langle \A\,\T\,\Pi\,f,f\rangle-\eps\,\langle \A\,\T\,(1-\Pi)\,f,f\rangle+\eps\,\langle \T\,\A\,f,f\rangle+\eps\,\langle \L\,f,(\A+\A^*)\,f\rangle\,.
\end{eqnarray*}
The proof of Theorems~\ref{Thm:Main1} and \ref{Thm:Main2} entirely
relies on the following estimate with $\eta=2\,\eps/(1-\eps)$.\smallskip
\begin{proposition}\label{Lem:HTheorem} Under the assumptions of
  Theorem~\ref{Thm:Main1} or \ref{Thm:Main2}, for any $\eps>0$ small
  enough, there exists an explicit constant $\lambda=\lambda(\eps)>0$ such that
  $D(f-F)+\lambda\,H(f-F)\le 0$ and $\liminf_{\eps\to 0}\lambda(\eps)/\eps>0$.
\end{proposition}

\section{Proofs of Proposition~\ref{Lem:HTheorem}}\label{Sec:Proofs}

To simplify the computations, we replace $f$ by $f-F$. Therefore, from now on we assume that $0=\int_{\R^d\times\R^d}f\,dx\,dv=\langle f,F\rangle$. By definition of $\Pi$, $\int_{\R^d}(\Pi\,f-f)\,dv=0$. We have obviously 
$\langle Lf,f \rangle \le - \|(1-\Pi) \, f \|^2$, and using the
identity $\eps\,x\,y\leq \frac c2\,x^2+\frac{\eps^2}{2\,c}\,y^2$, we
get the estimates
\begin{eqnarray*}
  &&-\eps\,\langle \A\,\T\,(1-\Pi)\,f,f\rangle=-\eps\,\langle \A\,\T\,(1-\Pi)\,f,\Pi\,f\rangle\le \frac{c_2}2\,\nrm{\A\,\T\,(1-\Pi)\,f}^2+\frac{\eps^2}{2\,c_2}\nrm{\Pi\,f}^2,\\
  &&\eps\,\langle \T\,\A\,f,f\rangle=\eps\,\langle \T\,\A\,f,(1-\Pi)\,f\rangle\le \frac\varepsilon{2}\,\nrm{\T\,A\,f}^2+\frac\eps{2}\,\nrm{(1-\Pi)\,f}^2,\\
  &&\eps\,\langle (\A+\A^*)\,\L\,f,f\rangle \le \eps\, \| (1-\Pi)\,f\|^2 + \eps\,\| \A\,f\|^2,
\end{eqnarray*}
for any $c_2>0$. Here we have used the identities
$\A\,\T\,(1-\Pi)=\Pi\,\A\,\T\,(1-\Pi)$ and $\T\,\A=(1-\Pi)\,\T\,\A$, which are respectively consequences of the fact that the range of $\A$ is contained in $\Pi\,L^2(d\mu)$ and that $\T\,\A\,f=v\,F\,\nabla_x(\rho(\A f)/\rho_F)$. We
moreover observe that $\langle \hat \a \cdot \a\, \Pi\,f,f\rangle =
\frac 1d\int_{\R^d} \big|\nx \big( \textstyle{\frac{\rho(f)}\rhoF}
\big) \big|^2\mF\,dx$. Let $g=\A\,f$, $u=\rho(g)/\rhoF$. An elementary
computation shows that $-\nx\int_{\R^d}v\,f\,dv=\rhoF\,u-\frac
1d\,\nx\left(\mF\,\nx u\right)$, from which one can deduce, after a few steps that we shall omit here, that
$\|\A\,f\|^2=\int_{\R^d}|u|^2\,\rhoF\,dx$ and $\|\T\,\A\,f\|^2=\frac
1d\int_{\R^d}|\nx u|^2\,\mF\,dx$ are such that
\[\label{Ineq:Estim}
2\,\|\A\,f\|^2+\|\T\,\A\,f\|^2\le\|(1-\Pi)\,f\|^2.
\]
We have proved that
\[
D(f) \le -\left(1-\frac 52\,\eps\right)\, \|(1-\Pi)\,f\|^2
-\eps\,\langle
\A\,\T\,\Pi\,f,f\rangle+\frac{c_2}2 \,\|\A\,\T\,(1-\Pi)\,f\|^2+\frac{\eps^2}{2\,c_2}\,\|\Pi\,
f \|^2.
\]
To complete the proof of Proposition~\ref{Lem:HTheorem}, it remains
to estimate from above $-\langle \A \,\T \,\Pi \,f,f \rangle$ and
$\|\A \, \T \, (1-\Pi) \, f \|^2$. As for the second of these two terms, we actually estimate $\|(\A\,\T\,(1-\Pi))^*f\|^2$ as follows.
 Using $(\A\,\T\,(1-\Pi))^*f=(\hat\a\cdot\a\,(1-\Pi))^*g$
with $g=\left( 1 + \hat\a \cdot\a\, \Pi \right)^{-1}f$, we first observe that
\be{Eqn:u1} \rho(f)=\rhoF\,u-\frac 1d\,\nx\left(\mF\,\nx u\right) \ee
where $u=\rho(g)/\rhoF$. Let $q_F:=\int_{\R^d}|v_1|^4\,F\,dv$, $u_{ij}:=\partial^2u/\partial x_i\partial x_j$. After some elementary but tedious computations, we also get
\be{Eqn:Dual}
\|(\A\,\T\,(1-\Pi))^*f\|^2=\frac 13\,\sum_{i,\,j=1}^d\,\int_{\R^d}\left[\left((2\,\delta_{ij}+1)\,q_F-\frac{3\,\mF^2}{d^2\,\rhoF}\,\delta_{ij}\right)u_{ii}\,u_{jj}+2\,(1-\delta_{ij})\,q_F\,u_{ij}^2\right]\,dx
\ee

\noindent\emph{Case of Theorem~\ref{Thm:Main1}.\/} In the \emph{Maxwellian case,\/}
various simplifications occur. With $\rhoF=e^{-V}=\frac 1d\,\mF=q_F$,
\eqref{Eqn:u1} becomes \be{Eqn:u2}
\rho(f)=u\,e^{-V}-\nx\left(e^{-V}\,\nx u\right) \ee and it follows
from (H3) that
\[
\langle
\A\,\T\,\Pi\,f,f\rangle\geq\frac\Lambda{1+\Lambda}\,\|\Pi\,f\|^2.
\]
On the other hand, from the above computation,
\[
\|(\A\,\T\,(1-\Pi))^*f\|^2 \le 2\,
   \sum_{i,\,j=1}^d\,\int_{\R^d}\left|u_{ij}\right|^2\,e^{-V}dx\,.
\]
Let $\|u\|_0^2:=\int_{\R^d}|u|^2\,e^{-V}dx$ and $W:=|\nx V|$. By
multiplying \eqref{Eqn:u2} by $u$, we get after an integration by
parts that $\|u\|_0^2+\|\nx u\|_0^2\le\|\Pi\,f\|^2$. By expanding the
square in $|\nx(u\,e^{-V/2})|^2$, one can prove using (H3) and (H4) that the following improved Poincar\'e inequality holds, with
$\kappa=(1-\theta)/(2\,(2+\Lambda\,c_0))$: \be{Improved Poincare}
\kappa\,\|W\,u\|_0^2\le\|\nx u\|_0^2 \ee for any $u\in H^1(e^{-V}dx)$
such that $\int_{\R^d}u\,e^{-V}dx=0$.

Multiply \eqref{Eqn:u2} by $W^2\,u$ and integrate by parts. By~(H5),
we get \be{Eq:interm} \|W\,u\|_0^2+\|W\,\nx u\|_0^2-2\,c_1\,\big(\|\nx
u\|_0+\|W\,\nx u\|_0 \big)\cdot\|W\,u\|_0\le \frac\kappa
8\,\|W^2u\|_0^2+ \frac 2\kappa\,\|\Pi\,f\|^2 \,. \ee Applying
\eqref{Improved Poincare} to $W\,u-\int_{\R^d}W\,u\,e^{-V}dx$, we get
\[
\kappa\,\|W^2u\|_0^2\le\int_{\R^d}|\nx(W\,u)|^2\,e^{-V}dx+2\,
\kappa\int_{\R^d}W\,u\,e^{-V}dx\int_{\R^d}W^3\,u\,e^{-V}dx\,.
\]
On the one hand, by the Cauchy-Schwarz inequality,
$\int_{\R^d}W\,u\,e^{-V}dx\le\|W\|_0\,\|u\|_0 =:a$, and on the other
hand, $\int_{\R^d}W^3\,u\,e^{-V}dx\leq a\,\|W\|_0^2+\frac
1{4\,a}\,\|W^2\,u\|_0^2$, so that
$2\int_{\R^d}W\,u\,e^{-V}dx\int_{\R^d}W^3\,u\,e^{-V}dx$ can be bounded
by $\frac 12\,\|W^2\,u\|_0^2+2\,\|W\|_0^4\,\|u\|_0^2$. Notice that
$\|W\|_0$ is bounded by (H6). As for the other term of the r.h.s., we
can simply write that $\int_{\R^d}|\nx(W\,u)|^2\,e^{-V}dx$ is bounded
by $2\,\|W\,\nx u\|_0^2+4\,c_1^2\,(\|u\|_0^2+\|W\,u\|_0^2)$ using
(H5). Hence we have
\[
\kappa\,\|W^2\,u\|_0^2\le 4\,\|W\,\nx
u\|_0^2+8\,c_1^2\,(\|u\|_0^2+\|W\,u\|_0^2)+4\,\kappa\,\|W\|_0^4\,\|u\|_0^2\,.
\]
Combined with \eqref{Eq:interm}, this proves that, for some $c_3>0$,
\[
\|W\,\nx u\|_0\le c_3\,\|\Pi\,f\|\,.
\]
By multiplying \eqref {Eqn:u2} by $\Delta u$ and integrating by parts,
we get
\[
\nrm{\nx^2 u}_0^2-\big(\nrm{W\,\nx
  u}_0+\nrm{\Pi\,f}\big)\nrm{\nx^2 u}_0\le\nrm{W\,u}_0\,\nrm{\nx
  u}_0\,.
\]
Altogether, this proves that $\|(\A\,\T\,(1-\Pi))^*f\|^2\le
c_4\,\|f\|^2$ for some $c_4>0$ and, as a consequence, $\|(\A\,\T\,(1-\Pi))f\|^2\le c_4\,\|f\|^2$. Since
$(1-\Pi)^2=1-\Pi$, we finally obtain
\[
\|(\A\,\T\,(1-\Pi))\,f\|^2\le c_4\,\|\,(1-\Pi)\,f\|^2.
\]
Summarizing, with $\lambda_1=1-\frac 12\,c_2\,c_4-5\,\eps/2$ and
$\lambda_2= \frac{\Lambda\,\eps}{1+\Lambda}-\frac{\eps^2}{2\,c_2}$, we
have proved that
\[
D(f)\le-\lambda_1\,\|(1-\Pi)\,f\|^2-\lambda_2\,\|\Pi\,f\|^2.
\]
With $c_2=a\,\varepsilon$, $a>\frac 12\,(1+1/\Lambda)$ and $\varepsilon>0$ small enough, $\lambda_1$ and $\lambda_2$ are positive and the result holds with $\lambda=\min\{\lambda_1,\lambda_2\}$. The explicit expression of $c_4$ can easily be retraced in
the above computations.

\medskip\noindent\emph{Case of Theorem~\ref{Thm:Main2}.\/} In the \emph{fast diffusion case,\/} we only sketch the main steps of the proof.
For $p=0$, $1$, $2$, let $w_p^2:=\omega_0\,V^{p-q}$, where
$q=k+1-d/2$, $w_0^2:=\rhoF$ and $V(x)=\left(1+|x|^2\right)^\beta$.
Define $\nrm u_i^2=\int_{\R^d}|u|^2\,w_i^2\,dx$. Notice that $\rhoF\in
L^1(\R^d)$ means $\beta(d-2k-2)+d<0$. This is the case if
$\beta(d+2-2k)+d-4<0$ and $\beta\ge 1$. The proof in the Maxwellian
case can be adapted as follows. Eq. \eqref{Eqn:u1} can be rewritten as
\be{Eqn:u3} \rho(f)=w_0^2\,u-\frac{2}{2k-d}\,\nx\left(w_1^2\,\nx u\right) \ee
and (H3) is replaced by the Hardy-Poincar\'e inequality, see
\cite{BlBoDoGrVa}, $\nrm u_0^2\le \Lambda\,\nrm{\nx u}_1^2$ for some
$\Lambda >0$, under the condition $\int_{\R^d}u\,w_0^2\,dx=0$. This
holds true if $\beta\ge 1$. The fact that $\langle
\A\,\T\,\Pi\,f,f\rangle\geq\frac\Lambda{1+\Lambda}\,\|\Pi\,f\|^2$ then
follows. We also need the following Hardy-Poincar\'e inequality
\[
\int_{\R^d}V^{\alpha+1-q-\frac
  1\beta}\,|u|^2\,dx-\frac{\Big(\;\int_{\R^d}V^{\alpha+1-q-\frac
    1\beta}\,u\,dx\Big)^2}{\int_{\R^d}V^{\alpha+1-q-\frac
    1\beta}\,dx}\le\frac
1{4 \, (\beta_0-1)^2}\int_{\R^d}V^{\alpha+1-q}\,|\nx u|^2\,dx
\]
which is responsible for the condition $\beta<\beta_0(\delta)$,
$\delta>0$. Observe that \eqref{Eqn:u3} multiplied by $u$ gives, after
an integration by parts, $\|u\|_0^2+(q-1)^{-1}\,\|\nx
u\|_1^2\le\|\Pi\,f\|^2$. By multiplying \eqref{Eqn:u3} by $V^\alpha\,
u$ with $\alpha:=1-1/\beta$ or by $V\,\Delta u$ and
integrating by parts, we find directly that $\|\nx^2 u\|_2^2$ is
bounded by $\|\Pi\,f\|^2$. Computations which are quite similar to the
ones of the Maxwellian case then allow to conclude. More details will
be given in \cite{DMS-1}.


\smallskip\noindent\copyright\ 2009 by the authors. This paper may be reproduced, in its entirety, for non-commercial purposes.

\smallskip\noindent{\bf Acknowledgments.} Partially supported by the
French-Austrian Amadeus project no. 13785UA, the ANR project IFO, the
Austrian Science Fund (project no. W8) and the European network DEASE. 
The authors thank an anonymous referee for his valuable comments and suggestions.

\medskip\begin{flushright}{\sl\today}\end{flushright}
\end{document}